\documentstyle[11pt]{article}

\input amssym.def
\input amssym

\def \a{{{\frak a}}}
\def \al{\alpha}
\def \ar{{\alpha_r}}
\def \Ad{{\rm Ad}}
\def \b{{{\frak b}}}
\def \bs{\backslash}
\def \C{{\bf C}}
\def \det{{\rm det}}

\def \E{{\cal E}}

\def \g{{{\frak g}}}
\def \ga{\gamma}
\def \Ga{\Gamma}
\def \h{{{\frak h}}}
\def \hol{{\rm hol}}

\def \k{{{\frak k}}}
\def \la{\lambda}
\def \lap{\triangle}
\def \m{{{\frak m}}}
\def \n{{{\frak n}}}
\def \N{\bf N}
\def \O{{\cal O}}
\def \p{{{\frak p}}}
\def \ph{\varphi}
\def \prf{{\bf Proof: }}
\def \qed{$ $\hfill {$\Box$}

$ $

}
\def \Q{\bf Q}
\def \R{{\bf R}}
\def \ra{\rightarrow}

\def \t{{{\frak t}}}
\def \tr{{\rm tr}}
\def \vol{{\rm vol}}
\def \V{{\cal V}}

\def \Z{\bf Z}

{\it}{}
\newtheorem{lemma}{Lemma.}[section]{\it}{}
{\it}{}
\newtheorem{proposition}{Proposition}[section]{\it}{}
\newtheorem{theorem}{Theorem}[section]{\it}{}
\newcommand{\rez}[1]{\frac{1}{#1}}

\newcommand{\binom}[2]{\left( \begin{array}{c}#1\\#2\end{array}\right)}

\begin{document}

\title{Holomorphic torsion for Hermitian locally symmetric spaces}
\author{by\\ {}\\ Anton Deitmar\\ {}\\  {\small Math. Inst., Im Neuenheimer Feld 288, 69120 Heidelberg, GERMANY}\\ {\it anton@mathi.uni-heidelberg.de}}

\date{}
\maketitle

\pagestyle{myheadings}
\markright{HOLOMORPHIC TORSION...}

\tableofcontents

\newpage

\begin{center}
{\bf Introduction}
\end{center}

The analytic torsion of the de Rham complex of an odd-dimensional manifold, introduced by D. Ray and I. Singer in \cite{RS-RT} was expressed as a special value of a zeta function defined in local geometric terms by D. Fried \cite{Fr} in case that the manifold was odd dimensional hyperbolic, i.e. the universal covering is an odd dimensional noncompact symmetric space of rank one.
Similar zeta functions have been studied for general rank one spaces by R. Gangolli \cite{Gang} and M. Wakayama \cite{Wak}.

The first breakthrough to higher rank was achieved by H. Moscovici and R. Stanton in \cite{MS-2}, where they used a higher supersymmetry argument for the calculation of orbital integrals. This method was developed further in \cite{D-Hitors} to cover all odd-dimensional spaces of higher rank. In this paper we give a construction for the Hermitian case that expresses the holomorphic torsion of Ray and Singer \cite{RS-HT} as a special value of a similar zeta function. See also \cite{Fr2} for the rank one case.

We get additional terms which may be interpreted as "$L^2$ determinants" and we calculate a
special value, which turns out to be a quotient of ordinary holomorphic torsion and the holomorphic $L^2$ torsion (see \cite{L}).

We will shortly describe the results: Let $X$ denote an Hermitian symmetric
space of the noncompact type. Let $\Ga$ be an
uniform torsion free lattice in the group of orientation preserving isometries
$G$ of $X$ and let $X_\Ga$ denote the quotient $\Ga \backslash X$. Write
$(\varphi ,V_\ph)$ for a finite dimensional unitary representation of $\Ga$.
For simplicity we will assume $X$ to be irreducible and $\Ga$ to be nice
\cite{Bor} which means that for each $\ga \in \Ga$ no eigenvalue of the
adjoint representation $\Ad (\ga)$ is a root of unity. 
From \cite{Bor} we take that every arithmetic $\Ga$ has a nice subgroup of
finite index.

To every Cartan subgroup $H$ of splitrank one we define the set ${\cal
  E}_H(\Ga)$ of $\Ga$-conjugacy classes that lie in $G$-conjugacy classes
meeting $H$ and let ${\cal E}_H^p(\Ga)$ denote the subset of primitive
conjugacy classes.
Since $\Ga$ is the fundamental group of $X_\Ga$ there is a bijection between
the set of conjugacy classes of $\Ga$ and the set of free homotopy classes of
closed paths in $X_\Ga$. 
For every $\ga \in \Ga$ let $X_\ga$ denote the union of all closed geodesics
in the class $[\ga]$. By \cite{DKV} we know that $X_\ga$ is a submanifold of
$X_\Ga$. For any manifold $Y$ let $\chi_{_1}(Y)$ denote the first higher Euler
number of $Y$, so $\chi_{_1}(Y) = -\sum_{p=0}^{\dim(Y)}p(-1)^p \dim(H^p(Y))$.

Write $A$ for the split part of the Cartan $H$ and fix a parabolic $P=MAN$.
Write $\n$ for the complex Lie algebra of $N$. Take $[\ga]\in \E_H(\Ga)$ then
$\ga$ is in $G$ conjugate to an element $h_\ga$ of $H$ which can be chosen to
act expandingly on $\n$. This action can be identified with the expanding part
of the Poincar\'e map around the geodesic $\ga$.

Let $c(H)$ denote the number of restricted roots with respect to the split
component of $H$, then $c(H)=1$ or $2$. In the case $c(H)=1$ we consider the
generalized Selberg zeta function
$$
Z_{H,1,\ph}(s) = \prod_{[\ga] \in \E_H^p(\Ga)} \prod_{N\geq 0}
\det(1-e^{-sl_\ga}\ga | V_\ph \otimes S^N(\n))^{\chi_{_1}(X_\ga)}.
$$

We show that $Z_{H,1,\ph}$ admits meromorphic continuation to $\C$ and indicate how to compute the poles and zeroes.
Set $Z_{H,\ph}(s) = Z_{H,1,\ph}(s+d(H))$ for a constant $d(H)$ given 
explicitly in the text (sec. \ref{constants}).

If $c(H)=2$ we set for $l\geq 0$:
$$
Z_{H,l,\ph}(s) = \exp \left( - \sum_{[\ga]\in \E_H(\Ga)}
\frac{\chi_{_1}(X_\ga) L^M(\ga ,\wedge^l\n_-)\tr(\ph(\ga))}
     {\det(1-\ga^{-1}|\n)}
\frac{e^{-sl_\ga}}
     {\mu_\ga} \right).
$$
The monodromy factor $L^M(\ga ,\wedge^l \n_-)$ is for regular $\ga$ equal to
$$
\frac{\tr(\ga | \wedge^l \n_-)}
     {\det(1-\ga^{-1}|\p_{M,+})}.
$$
Here $\p_M$ is the neutral part of the Poincar\'e map around $\ga$, which in this situation inherits a
complex structure of which $\p_{M,+}$ denotes the holomorphic component.
Then we set
$$
Z_{H,\ph}(s) = \prod_{l=0}^{\dim(\n_-)} Z_{H,l,\ph}(s+d_l(H))^{(-1)^l},
$$
for certain constants $d_l(H)$ which are given explicitly in \ref{constants}. Each $Z_{H,l,\ph}$ admits meromorphic continuation and we give formulas for the poles and zeroes. We then define
$$
Z_\ph(s) = \prod_H Z_{H,\ph}(s+b_0(H))
$$
for certain constants $b_0(H)$ given in \ref{constants}.

Let $n_0$ denote the vanishing order at zero of $Z_\ph(s)$ then our main result states that there is an explicit constant $c(X_\Ga ,\ph)$, which is a positive rational number, such that
 the function $R_\ph (s) = Z_\ph (s)s^{-n_0}/c(X_\Ga ,\ph)$ has the special value:
$$
R_{\varphi}(0) = \frac{T_\hol (X_\Ga ,\varphi)}{T_\hol^{(2)}(X_\Ga)^{\dim\ph}},
$$
where $T_\hol (X_\Ga ,\varphi)$ denotes the holomorphic torsion of the flat bundle defined by $\h$ and $T_\hol^{(2)}(X_\Ga)$ denotes the holomorphic $L^2$ torsion of $X_\Ga$.

One also can prove a formula for the case that $\Ga$ is not necessarily nice 
but still torsion free, which looks a little bit more complicated.

The contents of this paper also generalizes to give zeta functions for the 
higher holomorphic torsions $T_p$ or, more generally, to zeta functions for 
the torsion with coefficients in an arbitrary homogeneous bundle $E$.
The $p$-th holomorphic torsion is then the special case $E=\wedge^pT_+^*$, 
where $T_+$ is the holomorphic tangent bundle.

{\it I thank Ulrich Bunke, Andreas Juhl and Martin Olbrich for their comments and many stimulating discussions.}

$$ $$

{\bf Notation}

We will write $\N , \Z ,\Q ,\R ,\C$ for the natural, integer, rational, real and complex numbers.

For any locally compact group $G$ we write $\hat{G}$ for the unitary dual, i.e. the set of equivalence classes of irreducible unitary representations.

Throughout we will use small german letters to denote Lie algebras of the corresponding Lie groups which will be denoted by capital roman letters. A subscript zero will indicate the Lie algebra over $\R$, otherwise it will be its complexification, So for example $G$ a Lie group, then we write $\g_0=Lie_\R(G)$ and $\g =\g_0\otimes_\R \C$.

A {\bf virtual vector space} $V$ will be the formal difference of two vector spaces, i.e.:
$$
V = V^+ - V^-.
$$
An endomorphism $A$ of a virtual space $V$ is a pair of endomorphisms $ A^\pm $ on $ V^+ $ and $ V^- $ resp. The {\bf trace} and {\bf determinant} of $A$ are then
$$
\rm{tr} (A) = \rm{tr} (A^+) - \rm{tr} (A^-),\ \ \ \ \ \rm{\rm{det}} (A) = \rm{\rm{det}} (A^+) / \rm{\rm{det}} (A^-).
$$

Every vector space $V$ with $ \Z $-grading will naturally be considered as a virtual vector space by
$$
V^+ = V_{even},\ \ \ \ \ \ \ \ \ \ V^- = V_{odd}.
$$

As an example consider the exterior algebra over the finite dimensional space $V$ then for any endomorphism $A$ of $V$ we have the well known formula
$$
\rm{\rm{det}}( 1 - A ) = \rm{tr} (A \mid \wedge^* V ).
$$

For an endomorphism $A$ of a finite dimensional vector space $V$ denote by
$V_0 = \cup_{n\in \N} \ker A^n$ the generalized kernel of $A$. Since $V_0$ is stable under $A$ it induces an endomorphism of $V' = V/V_0$, which is injective.
Define the {\bf essential determinant} of $A$ to be
$$
\det'(A) = \det(A').
$$

\section{Holomorphic torsion}

We first need to lay some groundwork concerning determinants. A positive, possibly unbounded operator $A$ on a Hilbert space $H$ will be called {\bf zeta-admissible}, if
\begin{itemize}
\item[-] there is a $p\geq 1$ such that $A^{-p}$ is of trace class and
\item[-] the trace of the heat operator $e^{-tA}$ admits an asymptotic expansion
$$
\theta_A(t) = \rm{tr} e^{-tA} \sim \sum_{k=0}^\infty c_{\alpha_k} t^{\alpha_k}
$$
as $t\rightarrow 0$ where $\alpha_k \in \R$, $\alpha_k \rightarrow +\infty$.
\end{itemize}

The reader should notice that, from the fact that $A^{-p}$ is of trace class it immediately follows that the Hilbert space $H$ is separable and has a basis of $A$-eigenvectors, each eigenvalue occurring with a finite multiplicity and the eigenvalues cumulating at most at infinity.

Under these circumstances one defines the {\bf zeta function} of $A$ as:
$$
\zeta_A (s) = \rm{tr} A^{-s}\ \ \ \rm{for}\ \rm{Re}\ s > p.
$$
Then for $\rm{Re} s >> 0$ and $\lambda >> 0$ we consider the Mellin transform of $\theta_A$:
$$
M(z,\lambda )= \int_0^\infty t^{z-1} e^{-\lambda t} \theta_A(t)dt.
$$
The asymptotic expansion shows that for fixed $\lambda >0$ the function $z\mapsto M(z,\lambda ) = \Ga (z) \zeta_{A+\la}(z)$ is holomorphic in $\C$ up to simple poles at $-\alpha_k -n$, $k,n\geq 0$ of residue
$$
\rm{res}_{z=-\alpha_k-n} M(z,\lambda ) = \la^nc_k.
$$
(When two of those coincide the residues add up.)
This shows that $\zeta_{A+\la}(z)$ is regular at z=0. Extending the case of finite dimension we define:
$$
\rm{det}(A+\lambda ) = \rm{exp}(-\zeta_{A+\la}'(0)).
$$
Using the asymptotic expansion one sees that the function $\lambda \mapsto \rm{det}(A+\la )$ extends to an entire function with zeroes given by the eigenvalues of -A, the multiplicity of a zero at $\la_0$ being the multiplicity of the eigenvalue $\la_0$.

We will extend the definition of the determinant $\rm{det}(A)$ to the situation where not necessarily $A$ but  $A'=A\mid_{(\rm{ker}\ A)^\perp}$ is zeta admissible, then we will define $\zeta_A = \zeta_{A'}$ and $\det(A) = \det(A')$.

Let $E=(E_0 \rightarrow E_1 \rightarrow \dots \rightarrow E_n)$ be an elliptic complex over a smooth manifold $M$.
Assume each $E_k$ is equipped with a Hermitian metric. Then we can form the Laplace operators $\lap_k$ as second order differential operators.
When considered as unbounded operators on the spaces $L^2(E_k)$ these are known to be zeta admissible.
Now define the {\bf torsion} of $E$ as
$$
\tau_1 (E) = \prod_{k=0}^n \det(\lap_k)^{k(-1)^{k+1}}.
$$
Note that this definition differs by an exponent 2 from the original one \cite{RS-RT}.

For a compact smooth Riemannian manifold $M$ and $E\rightarrow M$ a flat Hermitian vector bundle the complex of $E$-valued forms on $M$ satisfies the 
conditions above so that we can define the torsion $\tau (E)$ via the de Rham complex. Now assume further, $M$ is K\"{a}hlerian and $E$ holomorphic then we may also consider the torsion $T(E)$ of the Dolbeault complex $\partial\ :\ \Omega^{0,.}(M,E) \rightarrow \Omega^{0,.+1}(M,E)$. 
The holomorphic torsion was introduced in \cite{RS-HT}.

We now define {\bf $L^2$-torsion}. For the following see also \cite{L}.
Let $M$ denote a compact oriented smooth manifold, $\Ga$ its fundamental group and $\tilde{M}$ its universal covering.
Let $E=E_0\rightarrow \dots \rightarrow E_n$ denote an elliptic complex over $M$ and $\tilde{E}=\tilde{E}_0\rightarrow \dots \rightarrow \tilde{E}_n$ its pullback to $\tilde{M}$.
Assume all $E_k$ are equipped with Hermitian metrics.

Let $\tilde{\triangle}_p$ and $\triangle_p$ denote the corresponding Laplacians.
The ordinary torsion was defined via the trace of the complex powers $\triangle_p^s$.
The $L^2$-torsion will instead be defined by considering the complex powers of $\tilde{\triangle}_p$ and applying a different trace functional.
Write ${\cal F}$ for a fundamental domain of the $\Ga$-action on $\tilde{M}$ then as a $\Ga$-module we have
$$
L^2(\tilde{E}_p) \cong l^2(\Ga) \otimes L^2(\tilde{E}_p\mid_{\cal F}) \cong l^2(\Ga) \otimes L^2(E_p).
$$
The von Neumann algebra $VN(\Ga)$ generated by the right action of $\Ga$ on $l^2(\Ga)$has a canonical trace making it a type $\rm{II}_1$ von Neumann algebra if $\Ga$ is infinite \cite{GHJ}. This trace and the canonical trace on the space $B(L^2(E))$ of bounded linear operators on $L^2(E)$ define a trace $\rm{tr}_\Ga$ on $VN(\Ga) \otimes B(L^2(E))$ which makes it a type $\rm{II}_\infty$ von Neumann algebra. The corresponding dimension function is denoted $\dim_\Ga$.  Assume for example, a $\Ga$-invariant operator $T$ on $L^2(E)$ is given as integral operator with a smooth kernel $k_T$, then a computation shows
$$
\tr_\Ga(T) = \int_{\cal F} \tr(k_T(x,x))dx.
$$ 
The reader familiar with the Selberg trace formula will immediately recognize this as the "term of the identity".

It follows for the heat operator $e^{-t\tilde{\lap}_p}$ that
$$
\tr_\Ga e^{-t\tilde{\lap}_p} = \int_{\cal F}\tr <x\mid e^{-t\tilde{\lap}_p}\mid x> dx.
$$
From this we read off that $\tr_\Ga e^{-t\tilde{\lap}_p}$ satisfies the same small time asymptotics as $\tr e^{-t\lap_p}$.

Let $\tilde{\lap}_p' = \tilde{\lap}_p|_{\ker(\tilde{\lap}_p)^\perp}$. Unfortunately very little is known about large time asymptotics of $\tr_\Ga(e^{-t\tilde{\lap}_p'})$ (see \cite{LL}).
Let
$$
NS(\lap_p) = \sup \{ \alpha \in \R \mid \tr_\Ga e^{-t \tilde{\lap}_p'} = O(t^{-\alpha/2})\ {\rm as}\ t\rightarrow \infty \}
$$
denote the {\bf Novikov-Shubin invariant} of $\lap_p$ (\cite{GrSh}, \cite{LL}).

Then $NS(\lap_p)$ is always $\geq 0$. J. Lott showed in \cite{L} that the Novikov-Shubin
invariants of Laplacians are homotopy invariants of a manifold. J.Lott and W. L\"uck conjecture in \cite{LL} that the Novikov-Shubin invariants of Laplace operators are always positive rational or $\infty$.

Throughout we will {\bf assume} that the Novikov-Shubin invariant of $\lap_p$ is positive. We will consider the integral
$$
\zeta_{\lap_p}^1(s) = \rez{\Ga (s)} \int_0^1 t^{s-1} \tr_\Ga e^{-t\tilde{\lap}_p'}dt,
$$
which converges for $\Re (s) >>0$ and extends to a meromorphic function on the entire plane which is holomorphic at $s=0$, as is easily shown by using the small time asymptotics (\cite{BGV},Thm 2.30).

Further the integral
$$
\zeta_{\lap_p(s)}^2(s) = \rez{\Ga(s)}\int_1^\infty t^{s-1}\tr_\Ga e^{-t\tilde{\lap}_p'} dt
$$
converges for $\Re (s)<\rez{2}NS(\lap_p)$, so in this region we define the {\bf $L^2$-zeta function} of $\lap_p$ as
$$
\zeta_{\lap_p}^{(2)} (s) = \zeta_{\lap_p}^1(s) + \zeta_{\lap_p}^2(s).
$$
Assuming the Novikov-Shubin invariant of $\lap_p$ to be positive we define the {\bf $L^2$-determinant} of $\lap_p$ as
$$
{\det}^{(2)}(\lap_p) = \exp (-\frac{d}{d{s}} \mid_{s=0} \zeta_{\lap_p}^{(2)} (s)).
$$

Now let the $L^2$-torsion be defined by
$$
T^{(2)}(E) = \prod_{p=0}^n {\det}^{(2)}(\lap_p)^{p(-1)^{p+1}}.
$$
 Again let $M$ be a K\"ahler manifold and $E\ra M$ a flat Hermitian holomorphic vector bundle then we will write $T^{(2)}_{hol}(E)$ for the $L^2$-torsion of the Dolbeault complex $\Omega^{0,*}(M,e)$.

\section{The trace of the heat kernel}

Let $\bar{X}$ denote a compact locally symmetric space whose universal covering $X$ is Hermitian globally symmetric of the noncompact type.
Then $X=G/K$ where $G$ is the connected component of the group of isometries of $X$ and $K$ is a maximal compact subgroup of $G$. Then there is a torsion-free uniform lattice $\Ga$ in $G$ such that $\bar{X} = \Ga \backslash X = \Ga \backslash G / K$. Further, $\Ga \cong \pi_1(\bar{X})$, the fundamental group of $\bar{X}$.
We will therefore write $X_\Ga$ instead of $\bar{X}$.
It follows that $G$ is a semisimple Lie group without center that admits a compact Cartan subgroup $T\subset K$.
We denote the real Lie algebras of $G$, $K$ and $T$ by $\g_0, \k_0$  and $\t_0$ and their complexifications by $\g ,\k$ and $\t$.
We will denote the Killing form of $\g$ by B.
As well, we will write $B$ for the diagonal of the Killing form, so $B(X) = B(X,X)$. Denote by $\p_0$ the orthocomplement of $\k_0$ in $\g_0$ with respect to $B$ then via the differential of \rm{exp} the space $\p_0$ is isomorphic to the real tangent space of $X=G/K$ at the point $eK$.
Let $\Phi (\t ,\g)$ denote the system of roots of $(\t ,\g)$, let $\Phi_c (\t ,\g) = \Phi (\t ,\k)$ denote the subset of compact roots and $\Phi_{nc} = \Phi -\Phi_c$ the set of noncompact roots.
To any root $\alpha$ let $\g_\alpha$ denote the corresponding root space.
Fix an ordering $\Phi^+$ on $ \Phi = \Phi (\t ,\g)$ and let $\p_\pm = \bigoplus_{\alpha \in \Phi_{nc}^+}\g_{\pm \alpha}$.
Then the complexification $\p$ of $\p_0$ splits as $\p = \p_+ \oplus \p_-$ and the ordering can be chosen such that this decomposition corresponds via \rm{exp} to the decomposition of the complexified tangent space of $X$ into holomorphic and antiholomorphic part.

Let $\theta$ denote the Cartan involution on $\g_0$ and on $G$ corresponding to the choice of $K$. Extend $\theta$ linearly to $\g$. Let $H$ denote a $\theta$-stable Cartan subgroup of $G$ then $H=A B$ where $A$ is the connected split component and $B$ compact.
The use of the letter $B$ here will not cause any confusion.
The dimension of $A$ is called the split rank of H. Let $\a$ denote the complex Lie algebra of $A$. Then $\a$ is an abelian subspace of $\p=\p_+ \oplus \p_-$. Let $X\mapsto X^c$ denote the complex conjugation on $\g$ according to the real form $\g_0$. The next lemma shows that $\a$ lies skew to the decomposition $\p =\p_+ \oplus \p_-$.

\begin{lemma} Let $Pr_\pm$ denote the projections from $\p$ to $\p_\pm$ 
then we have $\dim Pr_+(\a)$ = {\rm dim} $Pr_-(\a)$ = {\rm dim} $\a$, or, what amounts to the same: $\a \cap \p_\pm = 0$.
\end{lemma}

Proof: $\a$ is stable under complex conjugation which interchanges $\p_+$ and $\p_-$.
So let $X\in \p_+$ be such that $X+X^c \in \a$. Since $\a$ is abelian, the assumption $i(X-X^c)\in \a$ would lead to $0=[X+X^c,i(X-X^c)] = 2i[X^c,X]$ and the latter only
vanishes for $X=0$. So $\a$ does not contain $X$ or $X^c$ and the claim follows. \qed

Now let $\a'$ denote the orthocomplement of $\a$ in $Pr_+(\a) \oplus Pr_-(\a)$. For later use we write ${\cal V} = \a \oplus \a' = {\cal V}_+ \oplus {\cal V}_-$,
where ${\cal V}_\pm = {\cal V}\cap \p_\pm = Pr_\pm (\a)$.

Let $n=2m$ denote the real dimension of $X$ and for $0\leq p,q \leq m$ let $\Omega^{p,q}(X)$ denote the space of smooth $(p,q)$-forms on $X$.
Fix a finite dimensional unitary representation $(\varphi ,V_\varphi )$ of the group $\Ga$.
Then $\ph$ defines an Hermitian flat holomorphic vector bundle $E_\ph = X \times_\Ga V_\ph$ over $X_\Ga$.
On the space $\Omega^{p,q}(X_\Ga ,E_\ph )$ of $E_\ph$-valued forms we have a Laplacian $\triangle_{p,q,\Ga}^\ph$
and we write $\triangle_{p,q,\Ga}$ for the
Laplacian on $\Omega^{p,q}(X_\Ga)$.
Then $\triangle_{p,q,\Ga}$ is the pushdown of the Laplacian $\triangle_{p,q}$ on $\Omega^{p,q}(X)$.

By \cite{BM} the heat operator $e^{-t\triangle_{p,q}}$ has a smooth kernel $h_t^{p,q}$ of rapid decay in
$$
( C^\infty (G) \otimes \rm{End} (\wedge^p\p_+ \otimes \wedge^q\p_-))^{K\times K}.
$$
Now fix p and set for $t>0$
$$
f_t^p = \sum_{q=0}^m q(-1)^{q+1} \rm{tr}\ h_t^{p,q},
$$
where tr means the trace in $\rm{End} (\wedge^p\p_+ \otimes \wedge^q\p_-)$.

\subsection{Pseudocuspforms}\label{pseudocuspform}
We want to compute the trace of $f_t^p$ on the principal series representations.
To this end let $H=A B$ be a $\theta$-stable Cartan subgroup with split part $A$ and compact part $B$.
Let $P$ denote a parabolic subgroup of $G$ with Langlands decomposition $P=MA N$.
Note that $M$ is of inner type (see \cite{HC-HA1}, Lemma 4.9). Let $(\xi ,W_\xi)$ denote an irreducible unitary representation of $M$, $e^\nu$ a quasicharacter of $A$ and set $\pi_{\xi ,\nu} = Ind_P^G (\xi \otimes e^{\nu + \rho_{_P}}\otimes 1)$,
where $\rho_{_P}$ is the half of the sum of the $P$-positive roots.

\begin{proposition} The trace of $f_t^p$ under 
$\pi_{\xi ,\nu}$ vanishes if $\rm{dim}\ \a > 1 $. 
If $\rm{if\ dim\ } \a =1$ it equals
$$
e^{t\pi_{\xi ,\nu}(C)} \sum_{q=0}^{\dim(\p_-)-1} (-1)^q\rm{dim}(W_\xi \otimes \wedge^p\p_+ \otimes \wedge^q(\a^\perp \cap \p_-))^{K\cap M}.
$$
\end{proposition}

Proof: As before let $\a'$ the span of all $X-X^c$, where $X+X^c\in \a,\ X\in \p_+$ then ${\cal V} = \a \oplus \a' = \V_+ \oplus \V_-$ where $\V_\pm = \V \cap \p_\pm$.
The group $K_M=K\cap M$ acts trivially on $\a$ so for $x\in K_M$ we have $X+X^c = {\rm Ad}(x)(X+X^c) = {\rm Ad}(x)X + {\rm Ad}(x)X^c$.
Since $K$ respects the decomposition $\p=\p_-\oplus \p_+$, we conclude that $K_M$ acts trivially on $\V$, hence on $\V_-$. Let $r=\rm{dim}\a=\rm{dim}\V_-$. As a $K_M$-module we have
$$
\begin{array}{cll}
\wedge^p\p_-    &       =       &       \displaystyle
                                        \sum_{a+b=q} \wedge^a \V_- \otimes
                                        \wedge^b \V_-^\perp\\
                &       =       &       \displaystyle
                                        \sum_{a+b=q} \binom{r}{a}
                                        \wedge^b \V_-^\perp,
\end{array}
$$
where $\V_-^\perp = \a^\perp \cap \p_-$. By definition we get
$$
\begin{array}{cll}
\rm{tr}\pi_{\xi ,\nu}(f_t^p)    & = &   \displaystyle
                                \rm{tr} \pi_{\xi ,\nu}(\sum_{q=0}^m q(-1)^{q+1}
                                        h_t^{p,q})\\
                & = &   \displaystyle
                        e^{t\pi_{\xi ,\nu}(C)} \sum_{q=0}^m q(-1)^{q+1}
                                \rm{dim}(V_{\pi_{\xi ,\nu}} \otimes
                                \wedge^p\p_+ \otimes \wedge^q\p_-)^K
\end{array} $$
By Frobenius reciprocity this equals
\begin{eqnarray*}
                & = &   \displaystyle
                        e^{t\pi_{\xi ,\nu}(C)} \sum_{q=0}^m q(-1)^{q+1}
                                \rm{dim}(W_{\pi_{\xi}} \otimes
                                \wedge^p\p_+ \otimes \wedge^q\p_-)^{K\cap M}\\
                & = &   \displaystyle
                        e^{t\pi_{\xi ,\nu}(C)} \sum_{q=0}^m \sum_{a=0}^q q(-1)^{q+1}
                                \binom{r}{q-a} \rm{dim}(W_{\pi_{\xi}} \otimes
                                \wedge^p\p_+ \otimes
                                \wedge^a\V_-^\perp)^{K\cap M}\\
                & = &   \displaystyle
                        e^{t\pi_{\xi ,\nu}(C)} \sum_{a=0}^m \sum_{q=a}^m q(-1)^{q+1}
                                \binom{r}{q-a} \rm{dim}(W_{\pi_{\xi}} \otimes
                                \wedge^p\p_+ \otimes
                                \wedge^a\V_-^\perp)^{K\cap M}
\end{eqnarray*}
By taking into account $a\leq m-r$ we get  
$$
\sum_{q=a}^m q(-1)^q \binom{r}{q-a} = (-1)^r \binom{a}{1-r}
$$
 and
the claim follows. \qed

In view of the preceding proposition fix a $\Theta$-stable Cartan subgroup $H=A B$ with $\rm{dim}(A) =1$ and a parabolic $P=MA N$. Fix a system of positive roots compatible with $P$ and let $\rho$ denote the half sum of positive roots. For $\xi \in \hat{M}$ let $\la_\xi \in \b^*$ denote the infinitesimal character of $\xi$. Recall that we have
$$
\pi_{\xi ,\nu}(C) =B(\nu)+B(\la_\xi)-B(\rho).
$$

Let $c=c(H)$ denote the number of positive roots in $\phi(\a ,\g)$. Since $\a$ is a split torus and there is a real root, it follows $c=1$ or $c=2$.

\begin{lemma}
If $c=1$ let $\g =\g_1 \oplus \g_2$ be a decomposition of $\g$ into ideals where $\g_1$ is the simple ideal of $\g$ that contains $\a$. Accordingly we have $G\cong G_1 \times G_2$. Let $M_1 = M\cap G_1$ then we have an isomorphism as $K_M$-modules:
$$
\p_\pm \cong \C \oplus \p_{G_2 ,\pm} \oplus \p_{M_1}.
$$

In the case $c=2$ the space $M/K_M$ inherits the complex structure from $G/K$, i.e. $\p_M = \p_{M,+}\oplus \p_{M,-}$ where $\p_{M,\pm}=\p_M \cap \p_\pm$ and we have a $K_M$-module isomorphism:
$$
\p_\pm \cong \C \oplus \p_{M,\pm} \oplus \n_\pm,
$$
where
$$
n_\pm = \bigoplus _{\begin{array}{c}\alpha \in \phi_{nI}(\h ,\g)\\ \pm \alpha |_\b >0\end{array}} \g_\alpha.
$$
\end{lemma}

\prf
Let $H$ be a generator of $\a_0$. Write $\alpha_r$ for the unique positive real root in $\phi(\h ,\g)$. Since  for any root $\alpha$ we have $2B(\alpha ,\ar)/B(\ar) \in \pm \{0,1,2,3\}$ the only possible roots in $\phi^+(\a ,\g)$ are $\ar /2, \ar ,3\ar /2$.
Consider an embedding $\g \hookrightarrow gl_n$ such that the Cartan involution becomes $\theta(X)=-X^t$ and $\a$ is mapped to the diagonal.
Since $[\g_\ar ,\theta(\g_\ar)]=\a$ it is easy to see that $3\ar /2$ does not occur.

Now write $H=Y+Y^c$ with $Y\in \p_+$.
According to the root space decomposition $\g = \a \oplus \k_M \oplus \p_M \oplus \n \oplus \theta(\n)$ we write $Y=Y_a+Y_k+Y_p+Y_n+Y_{\theta(n)}$.
Because of $\theta(Y)=-Y$ it follows $Y_k=0$ and $Y_{\theta(n)}=-\theta(Y_n)$. For arbitrary $Z\in \b$ we have $[Z,Y]=0$ since $[Z,H]=0$ and the projection $Pr_+$ is $\b$-equivariant.
Hence $0=[Z,Y]=[Z,Y_a] + [Z,Y_p] + [Z,Y_n - \theta(Y_n)]$. Since $[Z,Y_a]=0$ and $[Z,Y_p]\in \p_M$ it follows $[Z,Y_n -\theta(Y_n)]\in \p_M \cap(\n \oplus \theta(n))=0$. Therefore $Y_p=0$
and $Y_n\in \n_r =\g_{\ar}$. Now $k_M$ acts trivially on $H$ and thus on $Y$, so on $Y_n$, i.e. $[\k_M ,Y_n]=0$.
Further $\a\cap \p_+$ is trivial, so $Y_n\neq 0$, so it generates the root space $\g_\ar$, so $[\k_M ,\g_\ar]=0$.

{\bf Assume} there is a root $\alpha \neq \ar$ such that $\alpha |_\a = \ar |_\a$. Then $B(\alpha ,\ar)>0$ and hence $\beta = \alpha -\ar$ is a root. The root $\beta$ is imaginary. Suppose $\beta$ is compact, then $\g_\beta \subset \k_M$ and we have $[\g_\ar ,\g_\beta]=\g_\alpha$, which contradicts $[\k_M,\g_\ar]=0$. It follows that $\beta$ is noncompact and thus it follows
$$
\g_{\ar |_\a}= [\p_M ,\g_\ar ] \oplus \g_\ar,
$$
where $\g_{\ar |_\a} = \{ X\in \g | [H,X] = \ar (H)X\}$ is the rootspace for the restricted root $\ar |_\a$. Let $\g_1$ as above the simple ideal containing $\a$, then $\g_1$ also contains $\g_\ar$ and hence we may substitute $\p_{M_1}$ in the above. Assume the root space $\g_{\ar /2 |_\a}$ is nonzero, then $(\a \oplus \m \oplus \n_r \oplus \g_{\ar /2|_\a} \oplus \g_{-\ar /2|_\a}) \cap \g_1$ would be a nontrivial ideal of $\g_1$, therefore $\ar /2$ does not occur as root of $(\a ,\g)$ and we may state
$$
\n = [\p_{M_1} ,\g_\ar] \oplus \g_\ar.
$$
The map $\phi : X \mapsto [X,Y_n]$ has no kernel in $\p_{M_1}$ since this would similarly allow us to construct a nontrivial ideal of $\g_1$. The map
\begin{eqnarray*}
\psi : \a \oplus \p_M \oplus \n & \ra  & \p \\
         a+p+n & \mapsto & a+p+n-\theta(n)
\end{eqnarray*}
is a $K_M$-isomorphism as is the map $\phi$. So, as $K_M$-modules
$\p \cong \a \oplus \n_r \oplus \p_M \oplus \p_M$. To derive the assertion for the case $c=1$ it remains to show $\p_{M_1}\cap \p_+ =0$, because then $\p_{M_1} \cong Pr_+(\p_{M_1})$. For this assume $X\in \p_{M_1}\cap \p_+$, then $0=[X,Y] = [x,Y_a] +[X,Y_n-\theta(Y_n)]$ and therefore $[X,Y_n]=0$ which implies $X=0$.

It remains the case where there exists no such root $\alpha$ as above.
This means $\n = \g_{\ar /2 |_\a} \oplus n_r$.
This gives $[\p_M ,\n_r]=0$ and thus $[\p_M,Y]=0$, which implies $\p_M =\p_{M,+} \oplus \p_{M,-} = \p_M \cap \p_+ \oplus \p_M \cap \p_-$. So $\p_M$ inherits the holomorphic structure and so does $\n$.
Since the maps $\phi$ and $\psi$ above are $K_M$-homomorphisms and $\p_+$ is given by a choice of positive roots the claim follows.
\qed

\subsection{Computation of Casimir eigenvalues}
In this section we will prove a general lemma which we will apply to the group $M$ later.
So here we will not assume $G$ connected, but $G$ should be of inner type (see \cite{HC-HA1} Lemma 4.9).
Also the space $X=G/K$ will not always be Hermitian.

Let $(\tau ,V_\tau)$ denote an irreducible unitary representation of $K$.
Since $G$ is of inner type, the Casimir operator $C_K$ of $K$ acts as a scalar $\tau(C_K)$ on $V_\tau$.

\begin{lemma} \label{pi(C)-lemma}
Assume $X$ Hermitian and let $(\pi ,W_\pi)$ be an irreducible unitary representation of $G$ and assume that
$$
\sum_{p=0}^{\dim \p_-} (-1)^p \dim(W_\pi \otimes \wedge^p\p_-\otimes V_\tau)^K \neq 0
$$
then we have
$$
\pi(C) = \tau(C_K) -B(\rho)+B(\rho_K).
$$
\end{lemma}

\prf
Consider $\p$ as a subspace of the Clifford algebra $Cl(B,\p)$. We will make $S=\wedge^*\p_-$ a $Cl(B,\p)$-module. For this let $x\in \p_-$ act on $S$ via
$$
x.z_1 \wedge \dots \wedge z_n = \sqrt{2}\ x\wedge  z_1 \wedge \dots \wedge z_n,
$$
and $y\in \p_+$ via
$$
y.z_1 \wedge \dots \wedge z_n = \sqrt{2} \sum_{i=1}^n (-1)^{i+1} B(y,z_i) z_1 \wedge \dots \hat{z_i}\dots \wedge z_n.
$$

This prescription turns $S$ into a nontrivial $Cl(B,\p)$-module. Since there is only one such of the dimension of $S$, we conclude that $S$ is the nontrivial irreducible $Cl(B,\p)$-module.
Therefore the argumentations of \cite{AtSch} apply to $S$ and especially the formula of Parthasarathy (\cite{AtSch}, (A13)).
Note that in \cite{AtSch} everything was done under the assumption of $G$ being connected. Since $G$ is of inner type, however, $\pi(C)$ will be a scalar. Writing $G^0$ for the connected component, the representation $\pi|_{G^0}$ will decompose as a finite sum of irreducibles on each of which the formula of Parthasarathy holds. Thus it holds globally.
Let $d_\pm$ be defined as in loc. cit. then our assumption leads to $ker(d_+ d_-) \cap \pi \otimes S(\tau) \neq 0$, and therefore $0=\tau(C_K)-\pi(C) -B(\rho) +B(\rho_K)$. \qed

\begin{lemma} \label{B(la_sigma)-B(rho)}
Assume $c=c(H)=2$ and let $(\sigma ,W_\sigma)$ be an irreducible $M$-subrepresentation of $\wedge^l\n$ such that $W_\sigma \cap \wedge^l\n_- \neq 0$.
Since $M$ is of inner type there is an infinitesimal character $\la_\sigma$ of $\sigma$ and we have
$$
B(\la_\sigma)-B(\rho)=-(1 +2\dim(\n_-)-l)^2B(\frac{\alpha_r}{2}),
$$
where $\alpha_r$ is the unique positive real root of $\t_M\oplus \a$.
\end{lemma}

\prf
The assertion follows from Corollary 5.7 of \cite{Kost}, here one uses the fact that the boundary map $\partial$ of the Lie algebra homology $H_*(\n)$ maps $\bigwedge^*(\n_+\oplus \n_-)$ to $\n_r {\wedge} \bigwedge^* (\n_+ \oplus \n_-)$ and maps $\bigwedge^*\n_-$ to zero. \qed

\begin{lemma}
Assume $c=c(H)=2$ and let $(\tau ,V_\tau)$ be an irreducible $K_M$-subrepresentation of $\wedge^l\n_-$, then the Casimir operator $C_{K_M}$ of $K_M$ acts on $V_\tau$ by the scalar
$$
\tau(C_{K_M}) = B(\rho_0) -  (1 +2\dim(\n_-)-l)^2B(\frac{\alpha_r}{2}) -\frac{2lB(\rho_{M_1,n})}{\dim(\p_{M_1,-})},
$$
where $\rho_{M_1,n}$ is the half of the sum of the noncompact roots of $\m_1=\m\cap \g_1$.
In the case that $rank(X)=1$, the last summand does not occur.
\end{lemma}

{\bf Convention.} In the case of rank one the last summand is undefined. To keep the formulas unified we will agree to consider it as zero then.

\prf
Let $x_1,\dots ,x_n$ be a basis of $\p_{M,+}$ consisting of root vectors with respect to $\t_M$ and such that $B(x_j ,x_j^c)=1$, where $x^c$ is the complex conjugate of $x$. Let $L$ be the element of the Lie algebra of $\k_M$ defined by
$$
L=\sum_j [x_j ,x_j^c].
$$

Since $(x_j^c)$ is the dual basis to $(x_j)$ it follows that $L$ lies in the center of $\k_M$. Furthermore $L$ is imaginary so that weights will take real values on $L$.

Consider the $K_M$-module $\n_-$.
We will show now that the center of $\k_M$ acts on $n_-$ by a character. 
To prove this we may momentarily assume that $\g$ is simple. It then follows that $\m$ is simple and that the centers of $\k$ and $\k_M$ are one dimensional.

Recall that $\t$ is spanned by $\t_M$ and $X+\theta(X)$, where $0\neq X\in \n_r$. Let $Z$ be a generator of the center of $\k$ then
$$
Z=a(X+\theta(X)) + bL +R
$$
where $R$ lies in $L^\perp \cap \t_M$. Since the latter is just the dual space of the span of the compact roots in $\t_M^*$ it follows $R=0$.

Now $Z$ acts as a scalar $\mu$ on $\psi(\n_-)$ and it is easy to see that $ab\neq 0$.
Let $Y\in \g_\alpha \subset \n_-$. A computation shows
$$
\mu Y = a[X,\theta(Y)] + b[L,Y].
$$
So it remains to show that $Y \mapsto [X,\theta(Y)]$ acts as a scalar on $\n_-$.
To this end let $H=[X,\theta(X)] \in \a$ then
$$
(\alpha_r(H)/2)Y = [H,Y] = [X,[\theta(X) ,Y]].
$$
Since root spaces are one dimensional there exists a number $c=c_Y$ such that $[X,\theta(Y)]=cY$. Therefore $[\theta(X),Y]=c\theta(Y)$ and $(\alpha_r(H)/2)Y=c^2Y$.
By this $c$ is determined up to sign. So $L$ acts on $\n_-$ by at most two eigenvalues.
Suppose there are two eigenvalues and an according decomposition
$$
\n_- = \n_-^0 \oplus \n_-^1,
$$
where $\n_-^0$ corresponds to the character of less absolute value.
We will show that $\g \cap ( \a \oplus \m \oplus \n_-^0 \oplus \n_-^{o\ c} \oplus \theta(\n_-^0) \oplus \theta(\n_-^{o\ c})$ is an ideal of $\g$ which cannot be since $\g$ is simple, thus proving that there is only one eigenvalue of $L$ on $\n_-$.
For this it suffices to show $[\p_{M,-} ,\n_-]=0$. Let $X\in \n_-$ then $X-\theta(X)$ is in $\p_-$ and therefore $[\p_{M,-},X-\theta(X)]=0$. But $\p_M$ also leaves stable the $\a$-rootspaces, therefore $[\p_{M,-},X]=0$.
We have shown that the center of $\k_M$ acts be a character $\chi_{_{\n_-}}$ on $\n_-$. As for the value of this character recall that the center of $\k_M$ acts by a character $\chi$ on $\p_{M,-}$.
The assumption $\chi \neq \chi_{_{\n_-}}$ would similarly allow us to construct a nontrivial ideal, hence it follows $\chi = \chi_{_{\n_-}}$.

Now we can drop the assumption of $\g$ being simple. 

Recall that each $x_j$ is a root vector, say $x_j \in \m_\alpha$ then $x_j^c\in \m_{-\alpha}$ and $[x_j,x_j^c]=H_\alpha$ where $H_\alpha \in \t_M$ is defined by $\alpha(H)=B(H,H_\alpha)$. Therefore
$$
L=\sum_{\alpha \in \phi_{noncompact}^+(\t_{M},\m)} H_\alpha.
$$

So let $\alpha $ be a noncompact negative root in $\phi(\t_{M_1},\m_1)$. Since $L$ acts on $\n_-$ by the same scalar as on $\p_{M_1,-}$ we get
$$
\chi_{\n_-}(L) = 2B(\rho_{M_1,n},\alpha) = -\frac{2B(\rho_{M_1,n})}{\dim(\p_{M_1,-})}.
$$

Now let $\mu$ be the lowest weight of $\tau$ and let $\sigma$ be the $M$-representation generated by $\tau$. We may assume that $\sigma$ is irreducible.
We claim that $\mu$ is the lowest weight vector of $\sigma$.
To see this it suffices to see that $\sigma(\p_{M,-})V_\tau =0$, which in turn follows from $\sigma(\p_{M,-})\n_- = [\p_{M,-},\n_-]=0$.
We thus have shown that $\mu$ is the lowest weight vector of $\sigma$. So by Lemma \ref{B(la_sigma)-B(rho)} we get
$$
B(\mu -\rho_M)-B(\rho) = -(1 +2\dim(\n_-) -l)^2B(\frac{\alpha_r}{2}).
$$
On the other hand we know $\tau(C_{K_M})=B(\mu -\rho_{K_M})-B(\rho_{K_M}).$ We have
$$
B(\mu - \rho_M) = B(\mu - \rho_{K_M}-\rho_{M,n}).
$$
Note that $2\rho_{M,n}$ is the dual of $L$. Since $L$ is in the center of $\k_M$ it follows
$$
B(\rho_{K_M} , \rho_{M,n})=0.
$$
Therefore
\begin{eqnarray*}
B(\mu -\rho_M) &=& B(\mu -\rho_{K_M}) -2B(\mu ,\rho_{M,n}) +B(\rho_{M,n}) \\
        &=& B(\mu -\rho_{K_M}) -\mu(L) +B(\rho_{M,n})\\
        &=& B(\mu -\rho_{K_M}) + \frac{2B(\rho_{M_1,n})}{\dim(\p_{M_1,-})}+B(\rho_{M,n}).
\end{eqnarray*}

From this the claim follows. \qed

Define $\hat{f_t^0}_{,H}$ as in \cite{HC-S}. Letting $b^*\in B^*$ we 
therefore get in the case that $c(H)=2$:
$$
\hat{f_t^0}_{,H}(\nu ,b^*) = e^{tB(\nu)}
$$ $$
\times \sum_{l=0}^{\dim(\n_-)} \sum_{p=0}^{\dim \p_{M,-}}(-1)^{p+l}
e^{-t((1+2\dim(\n_-)-l)^2B(\frac{\alpha_r}{2})+\frac{2B(\rho_{M_1,n})}
{\dim(\p_{M_1,-})}+B(\rho_{M,n}))} 
$$ $$
\dim(V_{b^*}\otimes \wedge^p \p_{M,-}\otimes \wedge^l\n_-)^{K_M}.
$$

$ $

In the other case, $c(H)=1$ we get

\begin{lemma}
Assume $c(H)=1$ then for $\pi_{\xi ,\nu}$ as in section \ref{pseudocuspform}
we get
$$
\tr(\pi_{\xi ,\nu}(f_t^0)) = e^{tC_\nu}
(\sum_{j=0}^{\dim(\p_{M_1})}(-1)^j\dim(W_{\xi_1}\otimes \wedge^j \p_{M_1})^{K\cap M_1})
$$ $$
\times (\sum_{k=0}^{\dim(\p_{G_2,-})}(-1)^k\dim(W_{\xi_2}\otimes \wedge^k\p_{G_2,-})^{K\cap G_2})
$$
with $C_\nu = B(\nu)+B(\rho_{K\cap G_2})+B(\rho_{M_1})-B(\rho)$
\end{lemma}

\prf
We have to show that $\tr(\pi_{\xi ,\nu}(f_t^0)) \neq 0$ implies 
$B(\la_\xi) = B(\rho_{M_1} + \rho_{K\cap G_2})$. 
We have $M=M_1 \times G_2$ and thus $\xi = \xi_1 \otimes \xi_2$. Therefore
$$
\tr(\pi_{\xi ,\nu}(f_t^0)) = e^{t\pi_{\xi ,\nu}(C)} \sum_{q=0}^{\dim(\p_-)-1} (-1)^q 
$$ $$
\sum_{j+k=q} \dim(W_{\xi_1}\otimes \wedge^j 
\p_{M_1})^{K\cap M_1} \dim(W_{\xi_2}\otimes \wedge^k\p_{G_2,-})^{K\cap G_2}
$$ $$
        = e^{t\pi_{\xi ,\nu}(C)} 
(\sum_{j=0}^{\dim(\p_{M_1})}(-1)^j\dim(W_{\xi_1}\otimes \wedge^j \p_{M_1})^{K\cap M_1}) 
$$ $$
\times (\sum_{k+0}^{\dim(\p_{G_2,-})}(-1)^k\dim(W_{\xi_2}\otimes \wedge^k\p_{G_2,-})^{K\cap G_2})
$$
Now assume $\tr(\pi_{\xi ,\nu}(f_t^0)) \neq 0$ then Lemma \ref{pi(C)-lemma} implies $B(\la_{\xi_2})=B(\rho_{K\cap G_2})$ whereas Lemma 2.4 in \cite{MS-2} gives $B(\la_{\xi_1})=B(\rho_{M_1})$. \qed

Therefore in the case $c(H)=1$ we get
$$
\hat{f_t^0}_{,H}(\nu, b^*) = \hat{f_t^0}_{,H}(\nu, b_1^*+b_2^*) 
$$ $$
= e^{t(B(\nu)+B(\rho_{K\cap G_2})+B(\rho_{M_1})-B(\rho))}
(\sum_{j=0}^{\dim(\p_{M_1})}(-1)^j\dim(V_{b_1^*}\otimes \wedge^j \p_{M_1})^{K\cap M_1}) \times
$$ $$
 (\sum_{k=0}^{\dim(\p_{G_2,-})}(-1)^k\dim(V_{b_2^*}\otimes \wedge^k\p_{G_2,-})^{K\cap G_2}).
$$

\subsection{The heat trace} 

Since $\Ga$ is the fundamental group of $X_\Ga$, every conjugacy class $[\ga]$ in $\Ga$ defines a free homotopy class of closed paths in $X_\Ga$. It is known \cite{DKV} that the union $X_\ga$ of all closed geodesics which are homotopic to $[\ga]$ is a smooth submanifold of $X_\Ga$. Let $\chi_{_1}(X_\ga)$ denote the first higher Euler number of $X_\ga$ (see \cite{D-Hitors}), i.e.
$$
\chi_{_1}(X_\ga) = -\sum_{p=0}^{\dim(X_\ga)} p(-1)^p b_p(X_\ga),
$$
where $b_p(X_\ga)$ is the $p$-th Betti number of $X_\ga$.

Now let $\E_H(\Ga)$ denote the set of nontrivial $\Ga$-conjugacy classes, which are in $G$ conjugate to an element of H. For $l\geq 0$ define
$$
b_l(H) = (\frac{c(H)}{2} + \dim(\n)-1-l)|\frac{\alpha_r}{c(H)}|.
$$
If $c(H)=2$ let
$$ \label{constants}
d_l(H) = \sqrt{b_l(H)^2  +\frac{2B(\rho_{M_1,n})}{\dim(\p_{M_1,-})}+B(\rho_{M,n})}.
$$
In the case $c(H)=1$ we finally set
$$
d(H) = \sqrt{B(\rho)-B(\rho_{K\cap G_2(H)})-B(\rho_{M_1(H)})}.
$$

For a finite dimensional unitary representation $(\ph ,V_\ph)$ of $\Ga$ and $\breve{\ph}$ denoting the dual of $\ph$
consider the unitary representation of $G$ given by right translations on
$$
L^2(\Ga \bs G,\breve{\ph}) = \{ f:G\rightarrow V_{\breve{\ph}} | f(\ga x) = \breve{\ph}(\ga)f(x)\forall
\ga \in \Ga ,\  
\int_{\Ga \bs G} \parallel f\parallel^2 < \infty \}.
$$
It is known that this representation splits discretely
$$
L^2(\Ga \bs G) = \bigoplus_{\pi \in \hat{G}} N_{\Ga ,\ph}(\pi) \pi
$$
with finite multiplicities $N_{\Ga ,\ph}\in \Z_{\geq 0}$. The Selberg trace
formula states that for a smooth function $f$ on $G$ which is rapidly
decreasing it holds
$$
\sum_{\pi \in \hat{G}} N_{\Ga ,\ph} (\pi) \tr (\pi(f)) = \sum_{[\ga]}
\vol(\Ga_\ga \bs G_\ga) \tr \ph(\ga)\O_\ga(f),
$$
where the sum runs over all conjugacy classes in $\Ga$ and $\O_\ga$ is the
orbital integral
$$
\O_\ga(f) = \int_{G/G_\ga} f(x\ga x^{-1}) dx.
$$

We want to compute the orbital integrals of the function $f_t^0$. 
Let at first $h\in G$ be a nonelliptic regular element. Since the trace of $f_t^0$ vanishes on principal series representations which do not come from splitrank one Cartan subgroups, we see that $\O_h(f_t^0)=0$ unless $h\in H$, a splitrank one Cartan.
Write $H=AB$ and $P=MAN$ as before. We have $h=a_hb_h$ and we call $h$ {\bf split-regular} if the split part $a_h$ is regular in $A$ which means that the centralizer of $a_h$ in $G$ equals the centralizer of $A$ in $G$. We will only be interested in the orbital integrals of split-regular elements.
We will further assume $a_h$ to be in the positive Weyl chamber $A^+$. 
For a split-regular $h$ we have by \cite{HC-DS},p. 32 ff with the notation from there:
$$
\O_h(f_t^0) = \frac{\tilde{\omega}_hF_{f_t^0}^H(h)}
                {[G_h:G_h^0] c_h h^\rho \det'(1-h^{-1}|(\g /\h)^+)}.
$$

Now assume first $c(H)=2$ then
Lemma 4.3 in \cite{MS-1} shows
$$
F_{f_t^0}^H(h) = \frac{e^{-l_h^2/4t}}{\sqrt{4\pi t}}h^{\rho_M}
\det(1-h^{-1}|(\k_M/\b)^+)
\sum_{l=0}^{\dim(\n_-)}(-1)^l e^{-td_l(H)^2}\tr (b_h | \wedge^l \n_-).
$$

From \cite{D-Hitors} we take in the case that $\Ga$ is nice:
$$
\chi_{_1}(X_\ga) = \frac{|W(\g_\ga ,\h)| \prod_{\alpha \in \Phi_\ga^+}(\rho_\ga ,\alpha)}               {\la_\ga c_\ga [G_\ga :G_\ga^0]}
        \vol(\Ga_\ga \bs G_\ga).
$$

Now let $(\tau ,V_\tau)$ be a finite dimensional unitary representation of 
$K_M$ and and define for $b\in B$ the monodromy factor:
$$
L^M(b,\tau) = \frac{\tilde{\omega}_b(b^{\rho_M} \det(1-b^{-1}|(\k_M/|\b)^+) \tr(\tau(b)))}
                {\tilde{\omega}_b(b^{\rho_M} \det(1-b^{-1}|(\m/\b)^+))}.
$$

Note that the expression $\tilde{\omega}_b(b^{\rho_M} \det (1-b^{-1}|(\m /\b)^+)) $ equals
$$
|W(\m_b ,\b)|\prod_{\alpha \in \Phi_b^+(\m ,\b)}(\rho_b,\alpha) b^{\rho_M}
\det'(1-b^{-1}|(\m /\b)^+),
$$
and that for $\ga =a_\ga b_\ga$ split-regular we have 
$$|W(\m_b,\b)|\prod_{\alpha \in \Phi_b^+(\m ,\b)}(\rho_b,\alpha) 
= |W(\g_\alpha ,\h)|\prod_{\alpha \in \Phi_\ga^+}(\rho_\ga ,\alpha),
$$ 
so that writing $L^M(\ga ,\tau) = L^M(b_\ga
  ,\tau)$ we get
$$
L^M(\ga ,\tau) = \frac{\tilde{\omega}_\ga(\ga^{\rho_M} \det (1-\ga^{-1} |
  (\k_M/\b)^+) \tr(\tau(b_\ga)))}
                      {|W(\g_\ga ,\h)|\prod_{\alpha \in \Phi_\ga^+}(\rho_\ga
                        ,\alpha) \det'(1-\ga^{-1} | (\m /\b)^+)}.
$$

Note further that for the virtual representation $\wedge^* \p_{M,+}$ it holds
$$
L^M(b ,\tau \otimes \wedge^* \p_{M,+}) = \tr(\tau(b)).
$$

From the above it finally follows that in the case $c(H)=2$ we get that
$\vol(\Ga_\ga \bs G_\ga) \O_\ga (f_t^0)$ equals:
$$
\frac{\chi_{_1}(X_\ga)}
     {\mu_\ga \det(1-\ga^{-1} | \n)}
\frac{e^{-l_\ga^2/4t}}
     {\sqrt{4\pi t}}
\sum_{l=0}^{\dim(\n_-)} (-1)^l e^{-td_l(H)^2} L^M(\ga ,\wedge^l \n_-).
$$

For the case $c(H)=1$ Lemma 4.3 in \cite{MS-1} gives
$$
F_{f_t^0}^H(h) = \frac{e^{-l_\h^2/4t}}
                      {\sqrt{4\pi t}}
 e^{-t(B(\rho_{K_2}) +B(\rho_{M_1})-B(\rho))}
 h^{\rho_M} \det (1-h^{-1} | (\k_M / \b)^+ \oplus \p_{M_1 ,+}).
$$
We can split up the differential operator $\tilde{\omega}_h$ according to the
decomposition $\g = \g_1 \oplus \g_2$ and analogous to the above we get in the
case $c(H)=1$ that $\vol(\Ga_\ga \bs G_\ga) \O_\ga (f_t^0)$ equals:
$$ 
\frac{\chi_{_1}(X_\ga)}
     {\mu_\ga \det(1-\ga^{-1} | \n)}
\frac{e^{-l_\ga^2/4t}}
     {\sqrt{4\pi t}}
e^{-t(B(\rho_{K_2})+B(\rho_{M_1})-B(\rho))} L^{G_2}(\ga ,1).
$$

Fix a finite dimensional unitary representation 
$(\ph ,V_\ph)$ of $\Ga$, let $E_\ph$ denote the flat hermitian bundle given by the dual $\breve{\ph}$ of $\ph$, i.e. $E_\ph := \Ga \bs X\times V_{\breve{\ph}}$.

\begin{theorem} \label{higher_heat_trace}
Let $X_\Ga$ be a compact locally Hermitian space with
  fundamental group $\Ga$ and such that the universal covering is globally
  symmetric
without compact factors. 
Assume $\Ga$ is nice and write $\triangle_{p,q,\ph}$ for the Hodge Laplacian
on (p,q)-forms with values in a the flat Hermitian bundle $E_\ph$, then
\begin{eqnarray*} & \displaystyle
\Theta (t) = \sum_{q=0}^{\dim_\C X_\Ga} q(-1)^{q+1} \tr\ e^{-t\triangle_{0,q,\ph}}
\\ & \displaystyle
= \sum_{\begin{array}{c}H/conj.\\ c(H)=2\end{array}} \sum_{[\ga] \in \E_H(\Ga)}
        \frac{\chi_{_1}(X_\ga) \tr \ph(\ga)}
                {\mu_\ga \det(1-\ga^{-1}|\n)}
        \frac{e^{-l_\ga^2/4t}}{\sqrt{4\pi t}}
\\ & \displaystyle
        \sum_{l=0}^{\dim(\n_-)} (-1)^l e^{-td_l(H)^2}
        L^M(\ga ,\wedge^l\n_-)
\\ & \displaystyle
+ \sum_{\begin{array}{c}H/conj.\\ c(H)=1\end{array}} \sum_{[\ga] \in \E_H(\Ga)}
        \frac{\chi_{_1}(X_\ga)\tr(\ph(\ga))}
             {\mu_\ga \det(1-\ga^{-1}|\n)}
 \frac{e^{-l_\ga^2/4t}}{\sqrt{4\pi t}}
                e^{-b_0(H)l_\ga} e^{td(H)^2} L^{G_2}(\ga ,1)
\\ & \displaystyle
        + f^0_t(e)\ \rm{dim}\ph\ vol(X_\Ga).
\end{eqnarray*}
\end{theorem}

The reader should keep in mind that by its definition we have for 
the term of the identity:
$$
f^0_t(e)\ \rm{dim}\ph\ vol(X_\Ga) = \sum_{q=0}^{\dim_\C X}q(-1)^{q+1} \tr_\Ga(e^{-t\lap_{0,q,\ph}}),
$$
where $\tr_\Ga$ is the $\Ga$-trace. Further note that by the Plancherel theorem the Novikov-Shubin invariants of all operators $\lap_{0,q}$ are positive.

\section{The zeta functions of Selberg and Ruelle} 

In this section we let $X$ to be an arbitrary globally symmetric space of the noncompact type. Fix a $\theta$-stable Cartan subgroup $H=A B$ of splitrank 1 and a parabolic subgroup $P$ with Langlands decomposition $P=MA N$.

\subsection{The generalized Selberg zeta function} 

Besides the parabolic $P=MAN$ we also consider the opposite parabolic $\bar{P} =MA\bar{N}$. The Lie algebra of $\bar{N}$ is written $\bar{\n}$. Let $V$ denote a Harish-Chandra module of $G$ then we consider the Lie algebra homology $H_*(\bar{\n},V)$ and cohomology $H^*(\bar{\n},V)$. It is shown in \cite{HeSch} that these are Harish-Chandra modules of the group $MA$.

We will say that a discrete subgroup $\Ga \subset G$ is {\bf nice} if for every $\ga \in \Ga$ the adjoint $\Ad(\ga)$ has no roots of unity as eigenvalues. Every arithmetic $\Ga$ has a nice subgroup of finite index \cite{Bor}.
Let $H_1\in A^+$ be the unique element with $B(H_1)=1$.

We will denote by ${\cal E}_H(\Ga)$ the set of nontrivial $\Ga$-conjugacy classes $[\ga]$, which are such that $\ga$ is in $G$ conjugate to an element of H. Such an element will then be written $a_\ga b_\ga$ or $a_\ga m_\ga$. The element $\ga \neq 1$ will be called {\bf primitive} if $\sigma \in \Ga$ and $\sigma^n =\ga$ with $\n\in \N$ implies $n=1$. Every $\ga \neq 1$ is a power of a unique primitive element. Obviously primitivity is a property of conjugacy classes. Let $\E_H^p(\Ga)$ denote the subset of $\E_H(\Ga)$ consisting of all primitive classes.

From \cite{D-onsome} we take:

\begin{theorem}\label{genSelberg}
Let $\Ga$ be nice and $(\ph ,V_\ph)$ a finite dimensional unitary representation of $\Ga$. Choose a $\theta$-stable Cartan $H$ of splitrank one. For $\Re(s)>>0$ define the zeta function
$$
Z_{H,\sigma,\ph}(s) = \prod_{[\ga]\in {\cal E}_H^p(\Ga)} \prod_{N\geq 0} \det(1-e^{-sl_\ga}\ga | V_\ph \otimes W_\sigma \otimes S^N(\n))^{\chi_{_1}(X_\ga)},
$$
where $\ga$ acts on $V_\ph \otimes W_\sigma \otimes S^N(\n)$ via $\ph(\ga) \otimes \sigma(m_\ga) \otimes Ad^N((m_\ga a_\ga)^{-1})$. Then $Z_{H,\sigma,\ph}$ has a meromorphic continuation to the entire plane. 
The vanishing order of $Z_{H,\tau ,\ph}(s)$ at a point $s=\la (H_1)$, $\la \in \a^*$, is
$$
(-1)^{\dim \ \n} \sum_{\pi \in \hat{G}}N_{\Ga ,\ph}(\pi)
\sum_{p,q}(-1)^{p+q} \dim (H^q(\bar{\n},\pi^0)\otimes \wedge^p\p_M \otimes V_{\breve{\tau}})^{K_M}_\la,
$$
where $(.)_\la$ denotes the generalized $\la$-eigenspace.
\end{theorem}

\prf: \cite{D-onsome}\qed

Note that $H^q(\n ,\pi^))$ and $H^q(\bar{\n},\pi^0)$ are Harish-Chandra modules for the group $MA$. The Cartan involution $\theta$ induces an isomorphism
$$
H^q(\bar{\n},\pi^0) \cong H^q({\n},\pi^0)(\theta),
$$
where the latter is $H^q({\n},\pi^0)$ with $\theta$-twisted action, so $ma$ would operate as $\theta (ma)$. Therefore the order of $Z_{H,\tau ,\ph}(s)$ at $s=\la (H_1)$ can also be expressed as
$$
(-1)^{\dim(\n)} \sum_{\pi \in \hat{G}} N_{\Ga ,\ph}(\pi) \sum_{p,q} (-1)^{p+q} \dim (H^q(\n ,\pi^0)\otimes \wedge^p\p_M \otimes V_{\breve{\tau}})^{K_M}_{-\la}.
$$

\begin{proposition} \label{extformel}
Let $\sigma$ be a finite dimensional representation of $M$ then the order of
$Z_{H,\sigma , \ph}(s)$ at $s=\la (H_1)$ is
$$
(-1)^{\dim(N)} \sum_{\pi \in \hat{G}}N_{\Ga,\ph}(\pi) \sum_{q=0}^{\dim(\m\oplus\n /\k_M)}(-1)^q \dim(H^q(\m\oplus\n, K_M ,\pi \otimes W_{\breve{\sigma}})_{-\la}).
$$
This can also be expressed as
$$
(-1)^{\dim(N)} \sum_{\pi \in \hat{G}}N_{\Ga,\ph}(\pi) \sum_{q=0}^{\dim(\m\oplus\n /\k_M)}(-1)^q \dim({\rm Ext}_{(\m \oplus \n ,K_M)}^q (W_\sigma ,V_\pi)_{-\la}).
$$
\end{proposition}

\prf \cite{D-onsome}\qed

\subsection{The generalized Ruelle zeta function} 

\begin{theorem}
Let $\Ga$ be nice and choose a $\theta$-stable Cartan subgroup $H$ of splitrank one. For $\Re(s)>>0$ define the zeta function
$$
Z_{H,\ph}^R(s) = \prod_{[\ga]\in {\cal E}_H^p(\Ga)} \det(1-e^{-sl_\ga}\ph(\ga)),
$$
then $Z_{H,\ph}^R(s)$ extends to a meromorphic function on $\C$.
\end{theorem}

\prf
Consider first the case when there is only one positive root $\alpha_r$ in the root system $\Phi(\a ,\g)$, we normalize $B$ such that $|\alpha_r|=2$. Then $a$ acts on $\wedge^l\bar{\n}$ by $a^{-l\alpha_r}$. Let $\sigma_l$ denote the representation of $K_M$ on $\wedge^l\bar{\n}$ the we get
$$
Z_{H,\ph}^R(s) = \prod_{l=0}^{\dim N} Z_{H,\sigma_l,\ph}(s+2l)^{(-1)^l}
$$
and hence the claim.

If we have two positive roots, say $\alpha_r$ and $\frac{\alpha_r}{2}$, with $|\alpha_r|=2$ then we have $\bar{\n} = \bar{\n}_r \oplus \bar{\n}_I$ where $\bar{\n}_r$ has dimension one and $a$ acts on $\bar{\n}_r$ by $a^{-\alpha_r}$ and on $\bar{\n}_I$ by $a^{-\frac{\alpha_r}{2}}$.
Now it follows
$$
Z_{H,\ph}^R(s) = \prod_{l=0}^{\dim N-1} (\frac{Z_{H,\wedge^l\bar{\n}_I,\ph}(s+l)}{Z_{H,\wedge^l\bar{\n}_I\otimes \bar{\n}_r,\ph}(s+l+2)})^{(-1)^l}. 
$$ \qed

\section{The holomorphic torsion zeta function}

Now let $X$ be Hermitian again and let $(\tau,V_\tau)$ be an irreducible unitary representation of $K_M$.

\begin{theorem}
Let $\Ga$ be nice and $(\ph ,V_\ph)$ a finite dimensional unitary
representation of $\Ga$. 
Choose a $\theta$-stable Cartan $H$ of splitrank one with $c(H)=2$. 
For $\Re(s)>>0$ define the zeta function
$$
Z_{H,\tau,\ph}^0(s) = \exp \left( -\sum_{[\ga] \in \E_H(\ga)}
\frac{\chi_{_1}(X_\ga) \tr(\ph(\ga)) L^M(\ga ,\tau)}
     {\det (1-\ga^{-1}|\n)}
\frac{e^{-sl_\ga}}{\mu_\ga})\right).
$$
Then $Z_{H,\tau,\ph}^0$ has a meromorphic continuation to the entire plane. 
The vanishing order of $Z_{H,\tau,\ph}^0(s)$ at $s=\la (H_1)$, $\la \in \a^*$ is
$$
-\sum_{\pi \in \hat{G}}N_{\Ga ,\ph}(\pi)
\sum_{p,q}(-1)^{p+q} \dim (H^q(\bar{\n},\pi)\otimes \wedge^p\p_{M,-} \otimes V_{\breve{\tau}})_\la^{K_M},
$$
where $(.)_\la$ means the generalized $\al$-eigenspace.
\end{theorem}

\prf
The proof proceeds as the proof of Theorem \ref{genSelberg} with 
$g_{\tau}$ instead of $f_\sigma$. 
\qed

Extend the definition of $Z_{H,\tau,\ph}^0(s)$ to arbitrary virtual representations in the following way. Consider a finite dimensional virtual representation $\xi = \oplus_i a_i \tau_i$ with $a_i\in \Z$ and $\tau_i\in \hat{K_M}$. Then let $Z_{H,\xi,\ph}^0(s) = \prod_i Z_{H,\tau_i,\ph}^0(s)^{a_i}.$

Now for $c(H)=2$ let $Z_{H,l,\ph}(s) = Z_{H,\wedge^l\n_-,\ph}^0(s)$
and
$$
Z_{H,\ph}(s) = \prod_{l=0}^{\dim(\n_-)} Z_{H,l,\ph}(s+d_l(H))^{(-1)^l},
$$

In the case $c(H)=1$ let
$$
Z_{H,\ph}(s) = Z_{H,1,\ph}(s+d(H)).
$$
See \ref{constants} for the constants.

For a nicer presentation of the results
 we will now assume that $G$ is simple.

\begin{proposition} \label{detformel}
Assume $\Ga$ is nice, then for $\la >>0$ we have the identity
$$
\prod_{q=0}^{\dim_\C X} (\frac{\det (\lap_{0,q,\ph}+\la)}{{\det}^{(2)} (\lap_{0,q,\ph}+\la)})^{q(-1)^{q+1}}
            =
$$ $$
\prod_{\begin{array}{c}H/{\rm conj.}\\ \rm splitrank=1\\ c(H)=2\end{array}} \prod_{l=0}^{\dim \n_+}  \left( Z_{H,l,\ph} (b_0(H) + \sqrt{\la + d_l(H)^2}\right)^{(-1)^l}
$$ $$
\prod_{\begin{array}{c}H/{\rm conj.}\\ \rm splitrank=1\\ c(H)=1\end{array}}   Z_{H,1,\ph} (b_0(H) + \sqrt{\la + d(H)^2}).
$$
\end{proposition}

\prf
Consider Theorem \ref{higher_heat_trace}.
For any semipositive elliptic differential operator $D_\Ga$ the heat trace $\tr e^{-tD_\Ga}$ has the same asymptotics as $t\ra 0$ as the $L^2$-heat trace $\tr_\Ga e^{-tD}$.
Thus it follows that the function
$$
h(t):= \sum_{q=0}^{\dim_\C X_\Ga} q(-1)^{q+1} (\tr e^{-t\lap_{0,q,\ph}}-\tr_\Ga e^{-t\lap_{0,q,\ph}})
$$
is rapidly decreasing at $t=0$.
Therefore, for $\la >0$ the Mellin transform of $h(t)e^{-t\la}$ converges for any value of $s$ and gives an entire function.
Let
$$
\zeta_\la (s) := \rez{\Ga (s)} \int_0^\infty t^{s-1} h(t) e^{-t\la} dt.
$$
We get that
$$
\exp(-\zeta_\la'(0)) = \prod_{q=0}^{\dim_\C X} (\frac{\det (\lap_{0,q,\ph}+\la)}{{\det}^{(2)} (\lap_{0,q,\ph}+\la)})^{q(-1)^{q+1}}.
$$
On the other hand, Theorem \ref{higher_heat_trace} gives a second expression for $\zeta_\la(s)$.
In this second expression we are allowed to interchange integration and summation for $\la >>0$ since we already know the convergence of the Euler products giving the right hand side of our claim.
\qed

Let $n_0$ be the order at $\la =0$ of the left hand side of the last proposition. Then
$$
n_0 = \sum_{q=0}^{\dim_\C(X)} q(-1)^q(h_{0,q,\ph} -h_{0,q,\ph}^{(2)}),
$$
where $h_{0,q,\ph}$ is the $(0,q)$-th Hodge number of $X_\Ga$ with respect to $\ph$ and $h_{0,q,\ph}^{(2)}$ is the $L^2$-analogue. Conjecturally we have $h_{0,q,\ph}^{(2)} = h_{0,q,\ph}$, so $n_0=0$. For $H$ a theta stable splitrank 1 Cartan with $c(H)=2$ and $l\geq 0$ let
$$
n_{H,l} = {\rm ord}_{s=b_0(H)+d_l(H)} Z_{H,l}(s),
$$
so
$$
n_{H,l} = -\sum_{\pi \in \hat{G}}N_{\Ga ,\ph}(\pi)
\sum_{p,q}(-1)^{p+q} \dim (H^q(\bar{\n},\pi)\otimes \wedge^p\p_{M,-} \otimes \wedge^l\n_+)_\la^{K_M}.
$$
Further, for $c(H)=1$ let $n_H$ be the order of $Z_{H,1,\ph}(s)$ at $s=b_0(H)+d(H)$, so
$$
n_H = -\sum_{\pi \in \hat{G}}N_{\Ga ,\ph}(\pi)
\sum_{p,q}(-1)^{p+q} \dim (H^q(\bar{\n},\pi)\otimes \wedge^p\p_{M})_\la^{K_M}.
$$
 We then consider
$$
c(X_\Ga ,\ph)= (\prod_{H, c(H)=2}\prod_{l=0}^{\dim \n_+}(2d_l(H))^{n_{H,l}(-1)^l})
(\prod_{H, c(H)=1} (2d(H))^{n_H}).
$$
Note that this is a positive rational number.

We assemble the results of this section to

\begin{theorem}
The zeta function $Z_{H,\ph}$ extends to a meromorphic function on the entire plane. Let
$$
Z_\ph(s) = \prod_{\begin{array}{c}H/{\rm conj.}\\ {\rm splitrank} = 1\end{array}} Z_{H,\ph}(s+b_0(H)).
$$
Let $n_0$ be the order of $Z_\ph$ at zero then Proposition \ref{detformel} shows that
$$
n_0 = \sum_{q=0}^{\dim_\C(X)} q(-1)^q (h_{0,q}(X_\Ga) -h_{0,q}^{(2)}(X_\Ga)),
$$
where $h_{p,q}(X_\Ga)$ is the $(p,q)$-th Hodge number of $X_\Ga$ and $h_{p,q}^{(2)}(X_\Ga)$ is the $(p,q)$-th $L^2$-Hodge number of $X_\Ga$.
Let $R_\ph(s) = Z_\ph(s)s^{-n_0}/c(X_\Ga ,\ph)$ then
$$
R_\ph(0) = \frac{T_\hol(X_\Ga ,\ph)}{T_\hol^{(2)}(X_\Ga)^{\dim\ph}}. 
$$ \qed
\end{theorem}

\tiny


\begin{thebibliography}{XX}

\begin{scriptsize}

\bibitem{AtSch} 
\bf Atiyah, M.; Schmid, W.: 
\it A Geometric Construction of the Discrete Series for Semisimple Lie Groups. \rm Invent. Math. 42, 1-62 (1977).

\bibitem{BM}  
\bf Barbasch, D.; Moscovici, H.:  
\it $L^2$-Index and the Selberg Trace Formula.  
\rm J. Func. An. 53, 151-201 (1983).

\bibitem{BGV}  
\bf Berline, N.; Getzler, E.; Vergne, M.:  
\it Heat Kernels and Dirac Operators.   
\rm Grundlehren 298. Springer 1992.

\bibitem{Bor}  
\bf Borel, A.:  
\it Introduction aux groupes aritm\'etiques.   
\rm Hermann. Paris 1969.

\bibitem{BorWall} 
\bf Borel, A.; Wallach,N.:  
\it Continuous Cohomology, Discrete Groups, and Representations of Reductive Groups.  
\rm Ann. Math. Stud. 94, Princeton 1980.

\bibitem{BrTD}  
\bf Br\"{o}cker, T.; tom Dieck, T.:  
\it Representations of Compact Lie Groups.   
\rm Springer 1985.

\bibitem{D-Det}  
\bf Deitmar, A.:  
\it A Determinant Formula for the generalized Selberg Zeta Function.   
\rm to appear in Quarterly J. Math.

\bibitem{D-Hitors}  
\bf Deitmar, A.:  
\it Higher torsion zeta functions.   
\rm Adv. Math. 110, 109-128 (1995).

\bibitem{D-Prod}  
\bf Deitmar, A.:  
\it Product expansions for zeta functions attached to locally homogeneous spaces.   
\rm Duke Math. J. 82, 71-90 (1996).

\bibitem{D-onsome}
\bf Deitmar, A.:
\it On some zeta functions attached to locally symmetric spaces of higher rank.
\rm preprint dg-ga/9511006.

\bibitem{DKV}  
\bf Duistermaat, J.J.; Kolk, J.A:C.; Varadarajan, V.S.:  
\it Spectra of locally
symmetric manifolds of negative curvature.   
\rm Invent. Math. 52 (1979) 27-93.

\bibitem{Fr}  
\bf Fried, D.:  
\it Analytic torsion and closed geodesics on hyperbolic manifolds.   
\rm Invent. Math. 84, 523-540 (1986).

\bibitem{Fr2}  
\bf Fried, D.:  
\it Torsion and closed geodesics on complex hyperbolic manifolds.   
\rm Invent. Math. 91, 31-51 (1988).

\bibitem{Gang} 
\bf Gangolli, R.:  
\it Zeta Functions of Selberg's Type for Compact Space Forms of Symmetric Spaces of Rank One.   
\rm Illinois J. Math. 21, 1-41 (1977).

\bibitem{GHJ}  
\bf Goodman, F.M.; de la Harpe, P.; Jones, V.F.R.:  
\it Coxeter Graphs and Towers of Algebras.   
\rm Springer 1989.

\bibitem{GrSh}  
\bf Gromov, M.; Shubin, M.A.:  
\it Von Neumann Spectra near Zero.   
\rm GAFA 1, 375-404 (1991)

\bibitem{HC-DS}  
\bf Harish-Chandra:  
\it Discrete series for semisimple Lie groups II.   
\rm Acta Math. 116, 1-111 (1966)

\bibitem{HC-HA1} 
\bf Harish-Chandra: 
\it Harmonic analysis on real reductive groups I. The theory of the constant term.  
\rm J. Func. Anal. 19 (1975) 104-204.

\bibitem{HC-S} 
\bf Harish-Chandra:  
\it Supertempered distributions on real reductive groups.  
\rm Studies in Appl. Math., Adv. in Math., Supplementary Studies Series, Vol 8,Acad. Press 139-153 (1983).

\bibitem{Helg}  
\bf Helgason, S.:  
\it Differential Geometry, Lie Groups, and Symmetric Spaces.   
\rm Acad. Press 1978.

\bibitem{Herb}  
\bf Herb, R.:  
\it Fourier inversion and the Plancherel Theorem.  
\rm in: Noncommutative harmonic analysis and Lie groups, SLN 880, 197-210 (1980).

\bibitem{HeSch}  
\bf Hecht, H.; Schmid, W.:  
\it Characters, asymptotics and $\n$-homology of Harish-Chandra modules.   
\rm Acta Math. 151, 49-151 (1983).

\bibitem{Hirz}  
\bf Hirzebruch, F.:  
\it Automorphe Formen und  der Satz von Riemann-Roch.   
\rm Symp. Int. Top. Alg. Mexico. 345-360 (1956).

\bibitem{Ju}  
\bf Juhl, A.:  
\it Zeta-Funktionen, Index-Theorie und hyperbolische Dynamik.   
\rm Habilitationsschrift. Humboldt-Universit\"{a}t zu Berlin 1993.

\bibitem{Knapp}  
\bf Knapp, A.:  
\it Representation Theory of Semisimple Lie Groups.   
\rm Princeton University Press 1986.

\bibitem{Koe}  
\bf K\"{o}hler, K.:  
\it Holomorphic torsion on Hermitian symmetric spaces.   
\rm to appear in: J. reine u. angew. Math.

\bibitem{Kost}  
\bf Kostant, B.:  
\it Lie algebra cohomology and the generalized Borel-Weil Theorem.   
\rm Ann. Math. 74 no 2 329-387 (1961).

\bibitem{Kot}  
\bf Kottwitz, R.:  
\it Tamagawa Numbers.   
\rm Ann. Math. 127, 629-646 (1988).

\bibitem{Lab}
\bf Labesse, J.P.:
\it Pseudo-coefficients tr\`es cuspidaux et K-th\'eorie.
\rm Math. Ann. 291, 607-616 (1991).

\bibitem{L}  
\bf Lott, J.:  
\it Heat kernels on covering spaces and topological invariants.   
\rm J.Diff. Geom. 35, 471-510 (1992)

\bibitem{LL}  
\bf Lott, J.; L\"uck, W.:  
\it $L^2$-Topologicla Invariants of 3-Manifolds.   
\rm Invent. Math. 120, 15-60 (1995).

\bibitem{MS-1}  
\bf Moscovici,H.; Stanton,R.:  
\it Eta invariants of Dirac operators on locally symmetric manifolds.   
\rm Invent. Math. 95, 629-666 (1989).

\bibitem{MS-2}  
\bf Moscovici,H.; Stanton,R.: 
\it R-torsion and zeta functions for locally symmetric manifolds.   
\rm Invent. Math. 105, 185-216 (1991).

\bibitem{RS-HT}  
\bf Ray, D.; Singer, I.:  
\it Analytic torsion for complex manifolds.   
\rm Ann. Math. 98, 154-177 (1973).

\bibitem{RS-RT}  
\bf Ray, D.; Singer, I.:  
\it R-torsion and the Laplacian on Riemannian manifolds.   
\rm Adv. Math. 7, 145-220 (1971).

\bibitem{Sel} 
\bf Selberg,A:  
\it Harmonic Analysis and Discontinuous Groups in weakly symmetric Riemannian spaces with Applications to Dirichlet Series.   
\rm J. Indian. Math. Soc. 20, 47-87 (1956).

\bibitem{Sh}  
\bf Shubin, M.:  
\it Pseudodifferential operators and Spectral Theory.   
\rm Springer 1987.

\bibitem{Voros}  
\bf Voros, A.:  
\it Spectral Functions, Special Functions and the Selberg Zeta Function.   
\rm Comm. Math. Phys. 110, 439-465 (1987).

\bibitem{Wak}  
\bf Wakayama, M.:  
\it Zeta functions of Selberg's type associated with homogeneous vector bundles.   
\rm Hiroshima Math. J. 15, 235-295 (1985).

\bibitem{Wall}  
\bf Wallach, N.:  
\it On the Selberg Trace Formula in the case of compact quotient.   
\rm Bull. AMS 82 No 2, 171-195 (1976).

\bibitem{Zel}  
\bf Zelobenko, D.P.:  
\it Compact Lie Groups and Their Representations.   
\rm Transl. Math. Monographs 40. AMS (1973)

\end{scriptsize}

\end{thebibliography}
\end{document}